\documentclass[12pt]{amsart}

\usepackage{amsmath,amsthm,amsfonts,xcolor}
\usepackage{tikz}
\usetikzlibrary{arrows}
\usetikzlibrary{decorations.markings}
\usetikzlibrary{shapes.geometric}
\usetikzlibrary{shapes.symbols}

\tikzstyle heightone=[scale=.8,shift={(0,-.3)}]
\tikzstyle heightones=[scale=.9,xscale=.4,shift={(0,.1)}]
\tikzstyle heightoneonehalf=[scale=.9,shift={(0,-.2)}]
\tikzstyle heighttwo=[scale=1,shift={(0,-.4)}]
\tikzstyle heighttwos=[scale=.4,xscale=.6,shift={(0,-.1)}]
\tikzstyle heightthree=[scale=.7,shift={(0,-.9)}]
\tikzstyle heightthrees=[scale=.5,xscale=.7,shift={(0,-.2)}]

\tikzstyle vertex=[circle,draw,fill=black,inner sep=1pt]
\tikzstyle ciliation=[circle,draw=none,fill=red,inner sep=1pt,semitransparent]
\tikzstyle ciliatednode=[vertex,pin={[pin distance=1mm,pin edge={semitransparent,red},ciliation]#1:{}}]

\tikzstyle matrix=[black,thick,signal,signal from=south,signal to=north,signal pointer angle=135,draw=blue!50,top color=blue!20,bottom color=black!10,scale=.8,inner sep=2pt]
\tikzstyle small matrix=[matrix,scale=.8]
\tikzstyle matrixonedge=[small matrix,sloped,rotate=-90]
\tikzstyle matrixonedge2=[small matrix,sloped,rotate=90]
\tikzstyle vector=[black,thick,rectangle,draw=gray!50!yellow,top color=yellow!30,bottom color=black!10,scale=.8,inner sep=2pt]
\tikzstyle small vector=[vector,scale=.8]
\tikzstyle plain vector=[rectangle,draw=none,fill=white,scale=.7]

\tikzstyle basiclabel=[draw=none,fill=none,shape=rectangle,inner sep=2pt,scale=.8]
\tikzstyle leftlabel=[basiclabel,anchor=east]
\tikzstyle rightlabel=[basiclabel,anchor=west]
\tikzstyle bottomlabel=[basiclabel,anchor=north]
\tikzstyle toplabel=[basiclabel,anchor=south]

\tikzstyle trivalent=[very thick]

\tikzstyle arrowstyle=[blue,semitransparent,scale=2]

\tikzstyle directed=[postaction={decorate,decoration={markings,
    mark=at position .65 with {\arrow[arrowstyle]{stealth}}}}]
\tikzstyle reverse directed=[postaction={decorate,decoration={markings,
    mark=at position .65 with {\arrowreversed[arrowstyle]{stealth};}}}]

\tikzstyle with matrix=[postaction={decorate,decoration={markings,
    mark=at position .5 with {\node[matrix]{#1};}}}]
\tikzstyle with small matrix=[postaction={decorate,decoration={markings,
    mark=at position .5 with {\node[small matrix]{#1};}}}]

\tikzstyle directed matrix=[postaction={decorate,decoration={markings,
    mark=at position .85 with {\arrow[arrowstyle]{stealth}},
    mark=at position .35 with {\node[matrix]{#1};}}}]
\tikzstyle directed small matrix=[postaction={decorate,decoration={markings,
    mark=at position .85 with {\arrow[arrowstyle]{stealth}},
    mark=at position .35 with {\node[small matrix]{#1};}}}]
\tikzstyle reverse directed matrix=[postaction={decorate,decoration={markings,
    mark=at position .4 with {\arrowreversed[arrowstyle]{stealth};},
    mark=at position .65 with {\node[matrix]{#1};}}}]
\tikzstyle reverse directed small matrix=[postaction={decorate,decoration={markings,
    mark=at position .4 with {\arrowreversed[arrowstyle]{stealth};},
    mark=at position .65 with {\node[small matrix]{#1};}}}]

\tikzstyle dotdotdot=[decorate,decoration={markings,
    mark=at position .3 with{\node{.};},
    mark=at position .5 with {\node{.};},
    mark=at position .7 with {\node{.};}}]

\tikzstyle wavyup=[out=90,in=-90]
\tikzstyle wavydown=[out=-90,in=90]

\tikzstyle permutation=[rectangle,fill=black,draw=black]
\tikzstyle symmetrizer=[rectangle,fill=black,draw=black]
\tikzstyle antisymmetrizer=[rectangle,fill=gray!10,draw=black]

\tikzstyle{every picture}=[semithick,baseline=0pt,heightone,label distance=-1mm]

\title{Unshackling Linear Algebra from Linear Notation}

\def\bfx{\mathbf{x}}
\def\bfu{\mathbf{u}}
\def\bfv{\mathbf{v}}
\def\bfw{\mathbf{w}}

\def\bse{\mathbf{\hat e}}
\def\sgn{\mathrm{sgn}}
\def\tr{\mathrm{tr}}

\def\det{\mathrm{det}}

\newtheorem{prop}{Proposition}

\newtheorem*{thm*}{Theorem}
\newtheorem*{fact*}{Fact}

\newtheorem*{matrixrule}{Matrix Rule}
\newtheorem*{noderule}{Node Rule}

\newtheorem*{inferrule}{Inferred Summation Rule}
\newtheorem*{freerule}{Free Edge Rule}

\theoremstyle{definition}
\newtheorem*{defn*}{Definition}

\theoremstyle{remark}

\newtheorem{example}[prop]{Example}
\newtheorem*{example*}{Example}
\newtheorem{exercise}[prop]{Exercise}
\newtheorem*{soln}{Solution}

\begin{document}

\ifnum0=1
\begin{abstract}
This paper introduces trace diagrams to an undergraduate audience, assuming only knowledge of basic linear algebra and graph theory. In particular, we describe how trace diagrams can be used to answer some natural questions that arise in the undergraduate mathematics curriculum: Why is the cross product only defined in three dimensions? What is the relationship between the cross product and the determinant? What's so special about the trace and the determinant of the matrix? What are the other coefficients that appear in the characteristic equation? These questions are answered using a diagrammatic algebra that offers visual insights, and that is in many respects a truer representation of linear algebra than the standard notation that historically developed.
\end{abstract}
\fi
%

\begin{center}
\Large
Unshackling Linear Algebra from Linear Notation
\end{center}

\bigskip

\begin{flushright}
Elisha Peterson\\
Department of Mathematical Sciences\\
United States Military Academy\\
646 Swift Road, West Point, NY 10996-1905\\
\verb+elisha.peterson@usma.edu+
\end{flushright}

\vspace{1in}

Have you ever seen one of those movies where the hero unearths an artifact covered with mysterious symbols, and it takes a brilliant scientist to decipher their meaning? The best of Hollywood creativity goes into devising a plausible system for the symbols. In the movie \emph{Contact}, for instance, Jodie Foster discovers a datastream of symbols that turn out to be a 3-dimensional blueprint for a device to contact alien life. Hollywood's tacit (and reasonable) assumption is that the mathematics of a different civilization would look very different.

This paper provides an accessible introduction to \emph{trace diagrams}, a non-traditional notation for linear algebra that could plausibly have been developed by another civilization. Trace diagrams are a completely different way of looking at vectors and matrices. Vectors are represented by \emph{edges} in a diagram, and matrices are \emph{markings} along the edges of the diagram. Instead of representing dot products and cross products by $\bfu\cdot\bfv$ and $\bfu\times\bfv$, one \emph{glues together} the edges corresponding to those vectors. One example is shown in Figure 1.

\begin{figure}[htb]
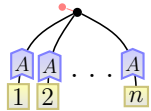
\label{f:td}
    $$
    \tikz[scale=.7,xscale=.7]{
        \node[ciliatednode=170](node)at(0,2){};
        \draw(-2,0)node[vector]{$1$}to[out=90,in=-160,looseness=1.2](node)node[pos=.25,small matrix]{$A$};
        \draw(-1,0)node[vector]{$2$}to[out=90,in=-135,looseness=1.2](node)node[pos=.3,small matrix]{$A$};
        \draw(-1,.5)to[dotdotdot](2,.5);
        \draw(2,0)node[vector]{$n$}to[out=90,in=-20,looseness=1.2](node)node[pos=.25,small matrix]{$A$};
    }
    $$
\caption{A trace diagram representation of a matrix determinant.}
\end{figure}

Surprisingly, the notation is perfectly rigorous, and often leads to proofs more elegant than those written using traditional notation, as we will show for vector identities. While trace diagrams have powerful applications in mathematics and physics \cite{Cvi08,Kau91}, the only prerequisite is an understanding of basic linear algebra and a willingness to work some examples to get used to doing real math with ``doodles''.

\subsubsection*{On Notation}
Whether one believes that mathematics is created or discovered, \textit{notation} is certainly created. And that notation can direct the course of mathematics. For instance, a matrix encodes the same information as a weighted directed graph, but the two notations inspire completely different questions. Matrices focus attention on concepts such as rank and invertibility, while graphs focus attention on nodes and flow between nodes.

Sometimes notations persist because they are useful and easy to understand, but there also cases where they exist simply because they are easier to write down in a single line of text. For example, permutation cycle notation such as $(2\:1\:3)$ is precise and easy to write down, but it is nonintuitive and often confuses beginners. The same permutation is much more clearly depicted by the diagram
    $$
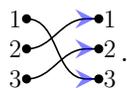
\tikz{
        \draw(0,0)node[vertex]{}node[leftlabel]{$3$}to[out=0,in=180](1,.5);\draw[directed](1,.5)to(1.2,.5)node[vertex,rightlabel]{$2$}node[vertex]{};
        \draw(0,.5)node[vertex]{}node[leftlabel]{$2$}to[out=0,in=180](1,1);\draw[directed](1,1)to(1.2,1)node[vertex,rightlabel]{$1$}node[vertex]{};
        \draw(0,1)node[vertex]{}node[leftlabel]{$1$}to[out=0,in=180](1,0);\draw[directed](1,0)to(1.2,0)node[vertex,rightlabel]{$3$}node[vertex]{};
    }.$$
There are many other cases, such as commutative diagrams, graph theory, and knot theory, where any attempt to fit the concepts into a single line would render them incomprehensible.

Outside of a few niche areas, diagrammatic notations have not been widely accepted. This is at least partly due to the difficulty in publishing. In fact, notations similar to trace diagrams have been used by mathematicians and physicists for years. In some cases researchers have ``invented'' new diagrammatic notations, only to discover that they almost perfectly replicate notation in a long-forgotten manuscript.

In the case of trace diagrams, there is a huge benefit in escaping the confines of a single line of text. The diagrammatic notation focuses attention on different questions, opening up new ideas that are not as readily achieved in the traditional manner. They also provide a refreshing perspective on the traditional theory.

\subsubsection*{Outline}

We begin with the definition of trace diagrams, and move directly into two special cases that help orient the reader to the diagrammatic point-of-view. We then provide an explicit description of how they are calculated. Finally, we provide the diagrammatic perspective on some questions often posed by students seeing vectors and linear algebra for the first time. We also look at some questions inspired by the diagrammatic notation.

We include several examples and exercises throughout, which are particularly important for adjusting to nonstandard notation. The reader is strongly encouraged to follow the examples closely and to complete the exercises.

\ifnum0=1
\subsection{On Mathematical Notation}

Whether one believes that mathematics is created or discovered, \emph{notation} is certainly created. And notation can direct the course of mathematics. For instance, a matrix encodes the same information as a weighted directed graph, but the two notations inspire completely different questions. So the body of mathematical knowledge depends in part upon the notation used in its construction.

Historically, mathematical notation has usually been developed along with new concepts. Good notation increases one's capacity to understand and reason with complex ideas. For this reason, most equations in textbooks are limited to ten or so different symbols. Imagine what calculus would be like if we had to write out the limit definition instead every time derivatives were mentioned!

Mathematical notations persist for three primary reasons: they enable higher-level thought, they are easy to understand, and they are easy to write down. The last constraint sometimes artificially forces complicated mathematics into ``linear notation'', a left-to-right, one-line format.\footnote{Frege's \emph{Begrifftschreit} is a rigorous diagrammatic notation for logic that failed in large part because it was difficult to typeset. Other notations, such as Russell's, were abandoned because they were incomprehensible. The common system of logic that eventually developed was both easy to typeset and easy to understand.} Because of this, notations can be more difficult to understand than an underlying idea. For example, the permutation $(1\:2\:3)$ in cycle notation is more naturally depicted by the diagram
    \tikz{\draw[directed](0,0)to[out=0,in=180](1,1);\draw[directed](0,.5)to[out=0,in=180](1,0);\draw[directed](0,1)to[out=0,in=180](1,.5);}.
There are many other cases, such as commutative diagrams, graph theory, and knot theory, where any attempt to fit the concepts into a single line would render them incomprehensible. More recently, computers have both made non-traditional notations easier to communicate and reinforced the need for discrete, character-based notations that computers can readily understand.

Trace diagrams and similar notations have been frequently rediscovered over the years, in large part due to the historic difficulty in typesetting ``nonlinear notation''.

\subsection{On Trace Diagrams}

The term \emph{trace diagram} was first rigorously defined in \cite{Pet06}, although the idea overlaps significantly with Penrose's \emph{tensor diagrams} \cite{Pen??}. Special cases of the diagrams have been used for years in various applications. The earliest such diagrammatic notations were used for angular momentum calculations \cite{Lev56}, and mathematicians have more recently explored their connections with knot theory \cite{Kau91}. The most comprehensive works on diagrammatic notations \cite{Blinn,Cvi08,Ste90} include a few minor results, but have minimal results regarding diagrams marked by matrices. Recent work on trace diagrams includes proofs of the Cayley-Hamilton Theorem \cite{Pet09a} and of several determinant formulas \cite{MP09}.
\fi

\section{General Trace Diagrams}\label{s:tracediagrams}

\begin{defn*}
    Given an underlying vector space of dimension $n$, a \emph{trace diagram} is a graph whose edges represent vectors. Endpoints of edges may be labeled by vectors, may be left \emph{free}, or may be joined at ordered nodes of degree $n$. Matrix markings may occur along any edge.
\end{defn*}
By \emph{ordered node}, we mean that the edges adjacent to the node are ordered. \emph{Matrix markings} are directed nodes placed along an edge and labeled by an $n\times n$ matrix.

The free ends of the diagram, along with those labeled by vectors, are typically partitioned into a set of \emph{inputs} and a set of \emph{outputs}. Each trace diagram gives rise to a function with corresponding input and output vectors.

If a trace diagram consists of two disconnected pieces, the individual functions are multiplied together to obtain the function of the entire diagram. Sums of multiple diagrams are also permitted, with the corresponding functions added together.

The next two sections describe special cases of trace diagrams, and are intended to clarify these ideas and to help the reader grow comfortable with diagrammatic reasoning. We will describe functions corresponding to trace diagrams, but will not say why these are the \textit{right} functions until section \ref{s:computation}, where we answer the question of how to compute functions for arbitrary diagrams.

The following result encapsulates the power of the diagrams:
\begin{fact*}
    Topologically equivalent trace diagrams correspond to equal functions, provided (i) the ordering at nodes is preserved, and (ii) the ordering of inputs and outputs is consistent.
\end{fact*}
Although the diagrams are drawn in the plane, no proof is necessary of this result since trace diagrams are defined combinatorially. We refer the reader to \cite{MP09} for further details.

\section{3-Diagrams}\label{s:vec}

The two main operations on 3-dimensional vectors are the \emph{cross product} and the \emph{dot product}. In 3-diagrams, a special case of trace diagrams that has been used for many years \cite{Ste90}, the diagrammatic forms of these operations are
    \begin{equation}\label{eq:d-dotcross}
	\bfu\times\bfv \leftrightarrow
    \tikz{\node[vertex]at(0,.5){}edge(0,1)edge[bend right]node[vector,pos=1]{$\bfu$}(-.5,0)edge[bend left]node[vector,pos=1]{$\bfv$}(.5,0);}
	\quad\text{and}\quad
	\bfu\cdot\bfv \leftrightarrow
    \tikz{\draw(0,.2)node[vector]{$\bfu$}to[out=90,in=90,looseness=2](1,.2)node[vector]{$\bfv$};}.
    \end{equation}

A \emph{$3$-diagram} is a graph with nodes of degree 3 and edges representing 3-vectors. If an edge is labeled by a particular vector, indicating that it has known value, it is an \emph{input} to the diagram's function. These are usually drawn at the bottom of the diagram. If an edge is \emph{free}, meaning it has an unlabeled end, it is an \emph{output} of the function, and usually drawn at the top of the diagram.

Both diagrams in \eqref{eq:d-dotcross} have input vectors $\bfu$ and $\bfv$. In the cross product diagram, the third strand is the output, identified with the vector result $\bfu\times\bfv$. The dot product diagram has no outputs. In this case, we say that the \emph{value} of the diagram is the scalar result $\bfu\cdot\bfv$.

\begin{example}\label{ex:3vec4}
	Draw the identity
	 $$(\bfu\times\bfv)\cdot(\bfw\times\bfx)=(\bfu\cdot\bfw)(\bfv\cdot\bfx)-(\bfu\cdot\bfx)(\bfv\cdot\bfw).$$
\begin{soln}
Keeping the vector inputs in the same order, the diagram is:
    $$
    
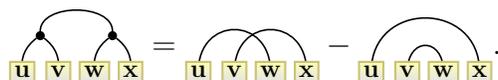
\tikz[scale=.6]{
        \node[vertex](node)at(.5,1){}
            edge[bend right]node[vector,pos=1]{$\bf u$}(0,0)
            edge[bend left]node[vector,pos=1]{$\bf v$}(1,0);
        \node[vertex](node2)at(2.5,1){}
            edge[bend right]node[vector,pos=1]{$\bf w$}(2,0)
            edge[bend left]node[vector,pos=1]{$\bf x$}(3,0)
            edge[bend right=90](node);
    }
    =
    \tikz[scale=.6]{
        \draw(0,0)node[vector]{$\bf u$}to[bend left=90,looseness=2](2,0)node[vector]{$\bf w$};
        \draw(1,0)node[vector]{$\bf v$}to[bend left=90,looseness=2](3,0)node[vector]{$\bf x$};
    }-
    \tikz[scale=.6]{
        \draw(0,0)node[vector]{$\bf u$}to[bend left=90,looseness=1.7](3,0)node[vector]{$\bf x$};
        \draw(1,0)node[vector]{$\bf v$}to[bend left=90,looseness=2.5](2,0)node[vector]{$\bf w$};
    }.
    $$
\end{soln}
\end{example}

\begin{exercise}\label{ex:diag-3vec}
What vector identity does this diagram represent?
    
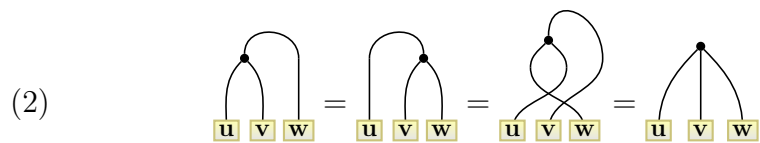
\begin{equation}\label{eq:diag-3vec}
    \tikz[scale=.6]{
        \node[vertex](node)at(.5,2){};
        \draw[](0,0)node[vector]{$\bf u$}to[out=90,in=-135](node);
        \draw[](1,0)node[vector]{$\bf v$}to[out=90,in=-45](node);
        \draw[](2,0)node[vector]{$\bf w$}to(2,2);\draw(2,2)to[out=90,in=90,looseness=1.5](node);
    }
    =\tikz[scale=.6]{
        \node[vertex](node)at(1.5,2){};
        \draw[](0,0)node[vector]{$\bf u$}to(0,2);\draw(0,2)to[out=90,in=90,looseness=1.5](node);
        \draw[](1,0)node[vector]{$\bf v$}to[out=90,in=-135](node);
        \draw[](2,0)node[vector]{$\bf w$}to[out=90,in=-45](node);
    }
    =\tikz[scale=.6]{
        \node[vertex](node)at(1,2.5){};
        \draw[](0,0)node[vector]{$\bf u$}to[out=90,in=-90](1.5,1.75);\draw(1.5,1.75)to[out=90,in=-45](node);
        \draw[](1,0)node[vector]{$\bf v$}to[out=90,in=-90](2.5,2)to(2.5,2)to[out=90,in=90,looseness=2](node);
        \draw(2,0)node[vector]{$\bf w$}to[out=90,in=-90](.5,1.75);\draw[](.5,1.75)to[out=90,in=-135](node);
    }
    =\tikz[scale=.7]{
        \node[vertex](node)at(1,2){};
        \draw[](0,0)node[vector]{$\bf u$}to[out=90,in=-135](node);
        \draw[](1,0)node[vector]{$\bf v$}to[out=90,in=-90](node);
        \draw[](2,0)node[vector]{$\bf w$}to[out=90,in=-45](node);
    }
    \end{equation}
(The reader is encouraged also to guess the meaning of the fourth term.)
\end{exercise}

This exercise illustrates the first kind of diagrammatic proof: the ``unproof''. No work is required to prove \eqref{eq:diag-3vec}, since the pictures are topologically equivalent and maintain the same ordering at the vertex. In contrast, here is the traditional direct proof of the first equality. Take $\bfu=(u_1, u_2, u_3)$ and a similar notation for $\bfv$ and $\bfw$. Then
\begin{align*}
    (\bfu\times\bfv)\cdot\bfw   &= (u_2v_3-u_3v_2, u_3v_1-u_1v_3, u_1v_2-u_2v_1)\cdot(w_1,w_2,w_3) \\
                                &= u_2v_3w_1 - u_3v_2w_1 + u_3v_1w_2 - u_1v_3w_2 + u_1v_2w_3 - u_2v_1w_3;\\
    \bfu\cdot(\bfv\times\bfw)   &= (u_1,u_2,u_3)\cdot(v_2w_3-v_3w_2, v_3w_1-v_1w_3, v_1w_2-v_2w_1) \\
                                &= u_1v_2w_3 - u_1v_3w_2 + u_2v_3w_1 - u_2v_1w_3 + u_3v_1w_2 - u_3v_2w_1 \\
                                &= (\bfu\times\bfv)\cdot\bfw.
\end{align*}
The proof is trivial, but completely unenlightening. A student might assume the identity is just a coincidence. In truth, the identity exists because of notation rather than the underlying mathematics.

A second type of proof is the ``surgery proof'', in which we use what we know about simple diagrams to perform manipulations with more complicated diagrams. For example, we know from Example \ref{ex:3vec4} that
    \begin{equation}\label{eq:basic3id1}
    \tikz[scale=.6]{
        \node[vertex](node)at(.5,1){}
            edge[bend right](0,0)
            edge[bend left](1,0);
        \node[vertex](node2)at(2.5,1){}
            edge[bend right](2,0)
            edge[bend left](3,0)
            edge[bend right=90](node);
    }
    =
    \tikz[scale=.6]{
        \draw(0,0)to[bend left=90,looseness=2](2,0);
        \draw(1,0)to[bend left=90,looseness=2](3,0);
    }-
    \tikz[scale=.6]{
        \draw(0,0)to[bend left=90,looseness=1.7](3,0);
        \draw(1,0)to[bend left=90,looseness=2.5](2,0);
    }.
    \end{equation}
Here, stating the identity without labelings indicates that the identity is true for any choice of vector labeling, under the assumption that inputs in the same position are labeled by the same vectors. This identity may be applied at \emph{any two adjacent nodes} in a 3-diagram.

\begin{example}
    Find an alternate expression of $(\bfu\times\bfv)\times\bfw$.
\begin{soln}
    First, represent $(\bfu\times\bfv)\times\bfw$ in diagrammatic form. Then apply \eqref{eq:basic3id1} in the neighborhood of the two nodes, manipulating the positions of the strands as appropriate:
    $$
    \tikz[scale=.6]{
        \node[vertex](node)at(.5,1.5){};\node[vertex](node2)at(1,2.5){};
        \draw[](0,0)node[vector]{$\bf u$}to[out=90,in=-135](node);
        \draw[](1,0)node[vector]{$\bf v$}to[out=90,in=-45](node);
        \draw(node)to[out=90,in=-135](node2);
        \draw[](2,0)node[vector]{$\bf w$}to[out=90,in=-45](node2);
        \draw(node2)to[out=90,in=-90](1,3.25);
    }
    =\tikz[scale=.6]{
        \draw[](0,0)node[vector]{$\bf u$}to[out=90,in=90,looseness=2](2,0)node[vector]{$\bf w$};
        \draw[](1,0)node[vector]{$\bf v$}to[out=90,in=-90](1,3.25);
    }
    -\tikz[scale=.6]{
        \draw[](0,0)node[vector]{$\bf u$}to[out=90,in=-90](1,3.25);
        \draw[](1,0)node[vector]{$\bf v$}to[out=90,in=90,looseness=2.5](2,0)node[vector]{$\bf w$};
    }.
    $$
    This proves the following identity:
    $$
    (\bfu\times\bfv)\times\bfw = (\bfu\cdot\bfw) \bfv - (\bfv\cdot\bfw) \bfu.
    $$
\end{soln}
\end{example}
By repeated application of \eqref{eq:basic3id1}, one can write \emph{any} single expression involving multiple cross products as a sum of expressions with at most one.

\section{Diagrams with Matrices}\label{s:diag-mx}

In diagrams with matrices, the \emph{matrix markings} are directed nodes placed along edges. We use a chevron symbol to represent matrices, as in the diagrams
	\begin{equation}\label{eq:mxmult}
    AB =
        \tikz[heightoneonehalf]{\draw(0,0)to(0,1.5)node[small matrix,pos=.5]{$AB$};}
        =\tikz[heightoneonehalf]{\draw(0,0)to(0,1.5)
            node[small matrix,pos=.33]{$B$}node[small matrix,pos=.67]{$A$};}
	\quad\text{and}\quad
	\bfv^T A\bfw =
        \tikz[heightoneonehalf]{
            \draw(0,0)node[vector]{$\bfw$}to[with small matrix={$A$}](0,1.5)node[vector]{$\bfv$};},
    \end{equation}
where $AB$ is a legal matrix product and $\bfv^T$ is a row vector. The directionality of the node is required to specify the input and output of the matrix.

One may use this to represent matrix entries diagrammatically, using the expression $a_{ij}=\bse^i A\bse_j$ for the entry in the $i$th row and $j$th column of $A$, where $\{\bse^i\}$ and $\{\bse_j\}$ are the standard row and column bases for the vector space. The diagram is
	\begin{equation}\label{eq:d-mxentry}
    a_{ij}=\bse^i A\bse_j =
        \tikz[heightoneonehalf]{\draw(0,0)node[vector]{$j$}to[with small matrix={$A$}](0,1.5)node[vector]{$i$};}.
    \end{equation}

\begin{example}\index{trace}
	Find a diagrammatic representation of the trace $\tr(A)$.
\begin{soln}
	$\tr(A) =  \sum_{i=1}^n a_{ii} = \sum_{i=1}^n
    \tikz[heightoneonehalf]{\draw(0,0)node[vector]{$i$}to[with small matrix={$A$}](0,1.5)node[vector]{$i$};}.$
\end{soln}
\end{example}
\begin{example}\index{determinant}
	Find a diagrammatic representation of the determinant
	\begin{equation}\label{detformula}
		\det(A)=\sum_{\sigma\in S_n}\sgn(\sigma)a_{1\sigma(1)}a_{2\sigma(2)}\cdots a_{n\sigma(n)}.
	\end{equation}
    Here, $\sigma\in S_n$ is a permutation on $n$ elements, and $\sgn(\sigma)$ is the \emph{signature} of the permutation, defined to be $(-1)^k$, where $k$ is the number of crossings in a diagram of the permutation.
\begin{soln}
One approach is to introduce new notation. If
    $\tikz[heightones,xscale=.6]{
        \foreach\xa in {1,2,5}{\draw(\xa,0)to(\xa,1);}
        \foreach\xa in {.1,.9}{\draw[dotdotdot](2,\xa)to(5,\xa);}
        \draw[antisymmetrizer](.7,.3)rectangle node[scale=.7]{$\sigma$}(5.3,.7);
    }$
represents a permutation $\sigma$ on the $n$ strands, then
    $$\det(A) = \sum_{\sigma\in S_n} \sgn(\sigma)
    \tikz[heighttwo,xscale=.5]{
        \foreach\xa/\xb in{1/1,2/2,5/n}{
            \draw(\xa,0)node[vector]{$\xb$}to(\xa,1)
                to[wavyup]node[small matrix,pos=.3]{$A$}(\xa,2)node[vector]{$\xb$};
        }
        \foreach\xa in {.25,1.75}{\draw[dotdotdot](2,\xa)to(5,\xa);}
        \draw[antisymmetrizer](.7,.5)rectangle node{$\sigma$}(5.3,0.9);
    }.
    $$
For example, if $n=2$, then
    $\det(A) =
    \tikz[heightoneonehalf,xscale=.5]{
        \foreach\xa/\xb in{1/1,2/2}{
            \draw(\xa,0)node[vector]{$\xa$}to[wavyup](\xb,.9)
                to[]node[small matrix,pos=.1]{$A$}(\xb,1.5)node[vector]{$\xb$};
        }}
    -
    \tikz[heightoneonehalf,xscale=.5]{
        \foreach\xa/\xb in{1/2,2/1}{
            \draw(\xa,0)node[vector]{$\xa$}to[wavyup](\xb,.9)
                to[]node[small matrix,pos=.1]{$A$}(\xb,1.5)node[vector]{$\xb$};
        }}
    = a_{11} a_{22} - a_{12} a_{21}.
    $
\end{soln}
\end{example}

\section{Computing Trace Diagram Functions}\label{s:computation}

Let $\{\bse_1,\ldots,\bse_n\}$ be an orthonormal basis for an $n$-dimensional vector space. A \emph{basis diagram} is an $n$-trace diagram in which every edge has been labeled by one of these basis elements. We use $i$ as shorthand for $\bse_i$.

The general process for computing a trace diagram's function involves two steps:
\begin{enumerate}
\item Express the diagram as a summation over basis diagrams.
\item Evaluate each basis diagram, and add the results together.
\end{enumerate}
We describe the second step first.

Basis diagrams evaluate to either 0 or a signed product of matrix and vector entries according to the following rules:
\begin{matrixrule}
    Each matrix and vector on the diagram contributes an entry to the product, according to the rules
        $$
        \tikz[heightoneonehalf]{\draw(0,0)node[vector]{$j$}to[with small matrix={$A$}](0,1.5)node[vector]{$i$};}=a_{ij}
        \qquad\text{and}\qquad
        \tikz{\draw(0,.1)node[vector]{$\bfu$}to(0,.9)node[vector]{$i$};}=u_i.
        $$
\end{matrixrule}

\begin{noderule}
    Nodes contribute multipliers of 0, if adjacent labels are repeated, or $\sgn(\tau)$, if the adjacent labels form a permutation $\tau$.
\end{noderule}

In the case where labels are not repeated, a convention is required to read off the permutation precisely. We use a small mark called a \emph{ciliation} near the node, and read off the permutation in a counter-clockwise fashion starting at that mark. (If the dimension is odd, the counter-clockwise ordering is enough to provide the sign of the permutation, since, for example, $\tbinom{1\:2\:3}{1\:2\:3}$ and $\tbinom{1\:2\:3}{2\:3\:1}$ have the same sign. So we can omit ciliations in this case.)

For example, the node
    $$
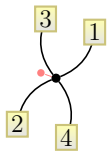
\tikz[shift={(0,.4)}]{
        \node[ciliatednode=170](va)at(0,0){}
            edge[bend right](50:1)edge[bend left](100:1)edge[bend right](230:1)edge[bend left](280:1);
        \node[vector]at(50:1){$1$};\node[vector]at(100:1){$3$};
        \node[vector]at(230:1){$2$};\node[vector]at(280:1){$4$};
    }$$
has edges ordered $(2,4,1,3)$, so the permutation is $\tbinom{1\:2\:3\:4}{2\:4\:1\:3}$. The number of transpositions required to express this permutation is 3, so the node contributes a sign of $-1$.

Now we treat the question of how to express diagrams as sums of basis diagrams. For now, we assume that there are no free edges. The summation uses the following rule:
\begin{inferrule}
    Express each edge in the diagram that is not already labeled by a basis element as a sum over basis labels:
    \begin{equation}\label{eq:summationrule}
        \tikz{\draw(0,0)to[wavyup](.1,1);} = \sum_{i=1}^n \tikz{\draw(0,-.15)to(.05,.25)node[vector,scale=.7]{$i$};\draw(.05,.75)node[vector,scale=.7]{$i$}to(.1,1.15);}.
    \end{equation}
\end{inferrule}
This rule is the analog of the expression of the identity matrix $I$ as the sum $I=\sum_{i=1}^n\bse_i\bse^i$.

The following example computes the dot product diagram in \eqref{eq:d-dotcross}.
\begin{example}
    Show that $\tikz{\draw(0,.2)node[vector]{$\bfu$}to[out=90,in=90,looseness=2](1,.2)node[vector]{$\bfv$};} = \bfu\cdot\bfv.$
\begin{soln}
    Use inferred summation over the interior edge.
    $$\tikz{\draw(0,.2)node[vector]{$\bfu$}to[out=90,in=90,looseness=2](1,.2)node[vector]{$\bfv$};}
    = \sum_{i=1}^n
    \tikz{\draw(.2,.1)node[vector]{$\bfu$}to(.2,.9)node[vector]{$i$};
          \draw(.8,.1)node[vector]{$\bfv$}to(.8,.9)node[vector]{$i$};}
    = \sum_{i=1}^n u_i v_i = \bfu\cdot\bfv.$$
\end{soln}
\end{example}

The next example justifies the terminology `trace' diagrams:
\begin{example}
    Evaluate the diagram $\tikz{\draw(0,.5)circle(.5);\node[small matrix]at(.5,.5){$A$};}$.
\begin{soln}
    Infer summation over the $n$ basic labels on the interior edge:
    \begin{equation}\label{eq:trace-diagram}
        \tikz{\draw(0,.5)circle(.5);\node[small matrix]at(.5,.5){$A$};}
        =\sum_{i=1}^n \tikz[heightoneonehalf]{\draw(0,0)node[vector]{$i$}to[with small matrix={$A$}](0,1.5)node[vector]{$i$};}
        =\sum_{i=1}^n a_{ii}
        =\tr(A).
    \end{equation}
\end{soln}
\end{example}

The rule
    $\tikz[heightoneonehalf]{\draw(0,0)to(0,1.5)
        node[small matrix,pos=.25]{$B$}node[small matrix,pos=.75]{$A$};}
    =\tikz[heightoneonehalf]{\draw(0,0)to(0,1.5)node[small matrix,pos=.5]{$AB$};}$
for matrix products follows by inferred summation on the strand connecting the two matrices:
    $$
    \tikz[heighttwo]{
        \draw(0,0)node[vector]{$j$}to(0,2)node[vector]{$i$}node[small matrix,pos=.33]{$B$}node[small matrix,pos=.67]{$A$};}
    = \sum_{k=1}^n
    \tikz[heighttwo]{
        \draw(0,0)node[vector]{$j$}to(0,1)node[vector]{$k$}node[small matrix,pos=.5]{$B$};
        \draw(-.6,1)node[vector]{$k$}to(-.6,2)node[vector]{$i$}node[small matrix,pos=.5]{$A$};}
    = \sum_{k=1}^n a_{ik} b_{kj}
    \equiv (AB)_{ij}
    = \tikz[heighttwo]{
        \draw(0,0)node[vector]{$j$}to(0,2)node[vector]{$i$}node[small matrix,pos=.5]{$AB$};}.
    $$
A similar approach allows vectors in diagrams to be reduced directly, in the sense that
    \begin{equation}\label{eq:d-vectorsumrule}
    \tikz[scale=.7]{\draw(1,0)node[vector]{$\bf v$}to(1,1);}
    =\sum_{i=1}^n v_i
    \tikz[scale=.7]{\draw(1,0)node[vector]{$i$}to(1,1);}.
    \end{equation}
This \emph{vector summation rule} is the analog of the equation $\bfv=\sum_{i=1}^n v_i \bse_i$, and may be proven by applying \eqref{eq:summationrule} on the edge in \tikz[scale=.7]{\draw(1,0)node[vector]{$\bf v$}to(1,1);}.

\begin{example}
    Given $\bfu=\tbinom{u_1}{u_2}$ and $\bfv=\tbinom{v_1}{v_2}$, evaluate
    $\tikz[scale=.7]{\node[ciliatednode=170](node)at(.5,1.5){};
            \draw(0,0)node[vector]{$\bf u$}to[out=90,in=-135](node);
            \draw(1,0)node[vector]{$\bf v$}to[out=90,in=-45](node);
        }.
    $
\begin{soln}
    Using \eqref{eq:d-vectorsumrule} and the node rule:
    $$\tikz[scale=.7]{\node[ciliatednode=170](node)at(.5,1.5){};
            \draw(0,0)node[vector]{$\bf u$}to[out=90,in=-135](node);
            \draw(1,0)node[vector]{$\bf v$}to[out=90,in=-45](node);
        }
        =\sum_{i=1}^2\sum_{j=1}^2 u_i v_j
        \tikz[scale=.7]{\node[ciliatednode=170](node)at(.5,1.5){};
            \draw(0,0)node[vector]{$i$}to[out=90,in=-135](node);
            \draw(1,0)node[vector]{$j$}to[out=90,in=-45](node);
        }
        = u_1 v_2
        \tikz[scale=.7]{\node[ciliatednode=170](node)at(.5,1.5){};
            \draw(0,0)node[vector]{$1$}to[out=90,in=-135](node);
            \draw(1,0)node[vector]{$2$}to[out=90,in=-45](node);
        }
        + u_2 v_1
        \tikz[scale=.7]{\node[ciliatednode=170](node)at(.5,1.5){};
            \draw(0,0)node[vector]{$2$}to[out=90,in=-135](node);
            \draw(1,0)node[vector]{$1$}to[out=90,in=-45](node);
        }
        = u_1 v_2 - u_2 v_1
        = \begin{vmatrix} u_1 & v_1 \\ u_2 & v_2 \end{vmatrix}.
    $$
\end{soln}
\end{example}

\begin{exercise}\label{ex:det3}
    Show, as suggested in Exercise \ref{ex:diag-3vec}, that
    $$
    \tikz[scale=.7]{
        \node[ciliatednode=170](node)at(1,2){};
        \draw[](0,0)node[vector]{$\bf u$}to[out=90,in=-135](node);
        \draw[](1,0)node[vector]{$\bf v$}to[out=90,in=-90](node);
        \draw[](2,0)node[vector]{$\bf w$}to[out=90,in=-45](node);
    }
    =
    \begin{vmatrix}
        u_1 & v_1 & w_1 \\
        u_2 & v_2 & w_2 \\
        u_3 & v_3 & w_3
    \end{vmatrix}
    =(\bfu\times\bfv)\cdot\bfw.
    $$
\end{exercise}

\begin{exercise}
    (Harder) Given an $n\times n$ matrix $A$, show that
    $$
    \tikz[scale=.7,xscale=.7]{
        \node[ciliatednode=170](node)at(0,2){};
        \draw(-2,0)to[out=90,in=-160,looseness=1.2](node)node[pos=.25,small matrix]{$A$};
        \draw(-1,0)to[out=90,in=-135,looseness=1.2](node)node[pos=.3,small matrix]{$A$};
        \draw(-1,.5)to[dotdotdot](2,.5);
        \draw(2,0)to[out=90,in=-20,looseness=1.2](node)node[pos=.25,small matrix]{$A$};
    }
    = \det(A)
    \tikz[scale=.7,xscale=.7]{
        \node[ciliatednode=170](node)at(0,2){};
        \draw(-2,0)to[out=90,in=-160](node);
        \draw(-1,0)to[out=90,in=-135](node);
        \draw(-1,.5)to[dotdotdot](2,.5);
        \draw(2,0)to[out=90,in=-20](node);
    }
    $$
    and that
    $$
    \tikz[heightoneonehalf]{
            \node[ciliatednode=140](topnode)at(0,1.5){};
            \node[ciliatednode=220](bottomnode)at(0,0){};
            \draw[with small matrix={$A$}](bottomnode)arc(270:90:.75);
            \draw[with small matrix={$A$}](bottomnode)arc(-90:90:.75);
            \draw[with small matrix={$A$}](bottomnode)to[out=135,in=-135](topnode);
            \draw[dotdotdot](-.25,.75)--(.65,.75);
        }
    = (-1)^{\lfloor\frac{n}{2}\rfloor} n! \det(A).
    $$
\end{exercise}

We have not yet described how to compute functions for diagrams with free (unlabeled) edges. These edges correspond to outputs of the underlying functions, and the rule for computing them is as follows:
\begin{freerule}
    A diagram with a single free edge corresponds to a function with output $\bfw=(w_1,w_2,\ldots,w_n)$, where $w_i$ is found by labeling the free strand by the basis label $i$ and evaluating the result.
\end{freerule}

For example,
    \tikz{\node[vertex]at(0,.5){}edge(0,1)edge[bend right]node[vector,pos=1]{$\bfu$}(-.5,0)edge[bend left]node[vector,pos=1]{$\bfv$}(.5,0);}
has a single free edge, corresponding to output $\bfw=(w_1,w_2,w_3)$. The vector summation and node rules give
    $$
    w_1=\tikz{
        \node[vertex]at(0,.5){}edge node[vector,pos=1]{$1$}(0,1)
            edge[bend right]node[vector,pos=1]{$\bfu$}(-.5,0)
            edge[bend left]node[vector,pos=1]{$\bfv$}(.5,0);}
    = u_2 v_3 - u_3 v_2.
    $$
\begin{exercise}
    Compute $w_2$ and $w_3$ in the above example. Then explain why
        $\tikz{\node[vertex]at(0,.5){}edge(0,1)edge[bend right]node[vector,pos=1]{$\bfu$}(-.5,0)edge[bend left]node[vector,pos=1]{$\bfv$}(.5,0);} = \bfu \times \bfv$.
\end{exercise}

The rules discussed in this section can be used to evaluate any trace diagram with a scalar or vector output. We refer the reader to \cite{MP09} for an explanation of how to compute diagrams with multiple outputs.

\ifnum0=1
    \footnote{For the interested reader, there are two basic approaches to proving rigor. The first technique first describes all the ``component pieces'' of diagrams, such as the value of a basic cup or cap, and then proves that combinations of these basic pieces are well-defined. The second technique uses \emph{signed graph coloring} to give a coefficient of a fully-labeled graph. The underlying functions are then described as multilinear functions whose coefficients are summations of these labeled graph coefficients.}
\fi

\section{A New Look at Some Old Questions}\label{s:questions}

The following questions are often asked by advanced students in multivariable calculus and linear algebra courses:
\begin{enumerate}
\item[Question 1:] Is there any pattern underlying vector identities with dot and cross products?
\item[Question 2:] Why is the cross product defined ``only'' in three dimensions?
\item[Question 3:] Is the rule for computing cross products using a determinant just a coincidence?
\item[Question 4:] Why do the trace and determinant show up in the discussion of eigenvalues and the $\det(A-\lambda I)=0$ equation?
\end{enumerate}
Some readers may have ready answers for these questions, but to a student seeing vectors or matrices for the first time, the answers may not be readily understood.

Trace diagrams are a good source of intuition for these problems, as it is often easier to generalize visual patterns than algebraic ones.

Question 1 has already been addressed. The diagrammatic rule \eqref{eq:basic3id1} unifies all vector identities involving more than one cross product. Finding the identity is as easy as finding two adjacent nodes in a diagram.

For Question 2, the visual cross product easily extends in dimension $n$ to
    $$
    \tikz{\node[ciliatednode=160]at(0,.5){}edge(0,1)
        edge[bend right]node[vector,pos=1,scale=.7]{$\bfu_1$}(-1,-.25)
        edge[bend right]node[vector,pos=1,scale=.7]{$\bfu_2$}(-.5,-.25)
        edge[bend left]node[vector,pos=1,scale=.7]{$\bfu_{n-1}$}(.85,-.25);
    \draw[dotdotdot](-.45,0)--(.75,0);
    }
    =\bfu_1\times\bfu_2\times\cdots\times\bfu_{n-1}
    .$$

From a visual perspective, there is zero coincidence in Question 3. Trace diagrams that differ only in placement of inputs and outputs represent the same underlying concept. In the correspondence with linear algebra, these inputs and outputs play almost identical roles in computation. So the determinant
    $
    \tikz{\node[vertex]at(0,.7){}edge[bend right]node[vector,pos=1]{$\bfu$}(-.5,.1)edge node[vector,pos=1]{$\bfv$}(0,.1)edge[bend left]node[vector,pos=1]{$\bfw$}(.5,.1);}
    =\det[\bfu\:\bfv\:\bfw]
    $
and the cross product
    $\tikz{\node[vertex]at(0,.5){}edge(0,1)edge[bend right]node[vector,pos=1]{$\bfu$}(-.5,0)edge[bend left]node[vector,pos=1]{$\bfv$}(.5,0);}
    =\bfu\times\bfv
    $
are merely two different windows into one idea. This visual perspective also begs the question of what the related diagrams
    $
    \tikz{\node[vertex]at(0,.5){}edge[bend left](-.5,1)edge[bend right](.5,1)edge node[vector,pos=1]{$\bfu$}(0,0);}
    $ and
    $
    \tikz{\node[vertex]at(0,.5){}edge[bend left](-.5,1)edge(0,1)edge[bend right](.5,1);}
    $
represent algebraically. As the answer involves terminology beyond the scope of this paper, we refer the interested reader to \cite{MP09} to uncover the answer

Before getting to Question 4, consider this combinatorial question: \emph{what are the simplest closed diagrams marked by a matrix $A$?} The simplest case is $\tikz{\draw[with small matrix={$A$}](0,.5)circle(.5);} = \tr(A)$, with no nodes. If the diagram has one node and is closed, $n$ must be even. Moreover, any unmarked edge in the diagram would produce a multiplier of 0, so the only case is
   \begin{equation}\label{eq:pfaffian-diagram}
   \tikz[heighttwo]{
        \node[ciliatednode=200]at(0,0){};
        \foreach\xa in{1.25,.6,.35}{\draw[with small matrix={$A$}](0,\xa)circle(\xa);}
        \draw[dotdotdot](-.5,1)to(-1.1,1.6);\draw[dotdotdot](.5,1)to(1.1,1.6);\draw[dotdotdot](0,1.2)to(0,2.5);
   }.
   \end{equation}
Two-node diagrams include the family
    \begin{equation}\label{eq:simplest-diagrams}
        \tikz[heightoneonehalf]{
            \node[ciliatednode=160](topnode)at(0,1.5){};
            \node[ciliatednode=200](bottomnode)at(0,0){};
            \draw[with small matrix={$A$}](bottomnode)arc(270:90:.75);
            \draw(bottomnode)arc(-90:90:.75);
            \draw(bottomnode)to[out=135,in=-135](topnode);
            \draw[dotdotdot](-.25,.75)--(.65,.75);
        }
        \qquad\cdots\qquad
        \tikz[heightoneonehalf]{
            \node[ciliatednode=160](topnode)at(0,1.5){};
            \node[ciliatednode=200](bottomnode)at(0,0){};
            \foreach\xa/\xb in{90/2.2,30/.5}{
                \draw[with small matrix={$A$},pos=.7](bottomnode)to[bend left=\xa,looseness=\xb](topnode);}
            \foreach\xa/\xb in{90/2,30/.5}{
                \draw(bottomnode)to[bend right=\xa,looseness=\xb](topnode);}
            \foreach\xa/\xb in{-.9/-.3,.2/.9}{\draw[dotdotdot](\xa,.75)to(\xb,.75);}
        }
        \qquad\cdots\qquad
        \tikz[heightoneonehalf]{
            \node[ciliatednode=160](topnode)at(0,1.5){};
            \node[ciliatednode=200](bottomnode)at(0,0){};
            \draw[with small matrix={$A$}](bottomnode)arc(270:90:.75);
            \draw[with small matrix={$A$}](bottomnode)arc(-90:90:.75);
            \draw[with small matrix={$A$}](bottomnode)to[out=135,in=-135](topnode);
            \draw[dotdotdot](-.25,.75)--(.65,.75);
        }.
    \end{equation}
The following result gives a partial answer to Question 4:
\begin{thm*}[Corollary 20 in \cite{Pet09c}]
    The diagrams in \eqref{eq:simplest-diagrams} are the coefficients of the characteristic polynomial $\det(A-\lambda I)$, up to a multiple that depends only on $n$ and the number of marked edges.
\end{thm*}
In this case, the diagrammatic pattern indicates a new way to understand the characteristic polynomial. The trace and determinant, which are constant multiples of the first and last diagrams in \eqref{eq:simplest-diagrams}, show up as special cases of a family of diagrams representing all coefficients of the polynomial. The traditional notation for these coefficients is $\tr_i(A)$, where $i=1,\ldots,n$ is the number of marked edges.

There is more to say about closed diagrams marked by a matrix $A$, and many unanswered questions. 
What is the function underlying \eqref{eq:pfaffian-diagram}? (We believe it to be the \emph{Pfaffian}, but have not yet proven this fact.) What are the other diagrams with two nodes, and what are their functions? What can be said about diagrams with more than two nodes?

\ifnum0=1
\section{Application to Trace Relations and Invariant Theory}\label{s:tracerelations}

Our final example leverages patterns in diagrammatic identities to produce trace relations for $2\times2$ matrices. All calculations are based upon the \emph{binor identity}
    \begin{equation}\label{eq:binor-mine}\index{binor identity}
    \tikz[xscale=.6]{\draw(0,0)to[wavyup](1,1);\draw(1,0)to[wavyup](0,1);}
    -\tikz[xscale=.6]{\draw(0,0)to[wavyup](.1,.5)to[wavyup](0,1);\draw(1,0)to[wavyup](.9,.5)to[wavyup](1,1);}
    -\tikz[xscale=.6]{\draw(0,0)to[out=90,in=90](1,0);\draw(0,1)to[out=-90,in=-90](1,1);
        \node[ciliatednode=170] at(.5,.3){};\node[ciliatednode=190] at(.5,.7){};}
    =0,
    \end{equation}
which is the two-dimensional analog of \eqref{eq:basic3id1}. We also need the fact that
    $\tikz{\draw(0,0)to(0,.3)node[ciliatednode=180]{}to(0,.7)node[ciliatednode=180]{}to(0,1);}
        = -\tikz{\draw(0,0)to(0,1);}$
    but $\tikz{\draw(0,0)to(0,.3)node[ciliatednode=180]{}to(0,.7)node[ciliatednode=0]{}to(0,1);}
        = \tikz{\draw(0,0)to(0,1);}$, and the identity
    \begin{equation}\label{eq:detnode}
    \tikz[xscale=.6]{\draw(-.5,.1)to(-.5,.5)to[out=90,in=90](.5,.5)to(.5,.1);\node[small matrix]at(-.5,.4){$A$};
            \node[small matrix]at(.5,.4){$A$};\node[ciliatednode=160]at(0,.78){};}
        = \det(A) \tikz[xscale=.6]{\draw(-.5,.1)to(-.5,.5)to[out=90,in=90](.5,.5)to(.5,.1);\node[ciliatednode=160]at(0,.78){};}.
    \end{equation}
\begin{exercise}
    Prove \eqref{eq:binor-mine} and \eqref{eq:detnode}
\end{exercise}

This is an extremely powerful result, since it can be used to remove all crossings from a 2-diagram. The characteristic equation pops out when we add in some matrix labels and connect the strands:
    \begin{equation}\label{eq:ch2a-diagram}
        \tikz[xscale=.6]{
            \draw(0,-.3)--(0,0)to[wavyup](1,1)node[small matrix]{$A$}to[out=90,in=90](2,1)--(2,0)to[out=-90,in=-90](1,0);
            \draw(1,0)to[wavyup](0,1)node[small matrix]{$A$}--(0,1.4);}
        -\tikz[xscale=.6]{
            \draw(0,-.3)--(0,0)to[wavyup](.1,.5)to[wavyup](0,1)node[small matrix]{$A$}--(0,1.4);
            \draw(1,0)to[wavyup](.9,.5)to[wavyup](1,1)node[small matrix]{$A$}to[out=90,in=90](2,1)--(2,0)to[out=-90,in=-90](1,0);
            }
        -\tikz[xscale=.6]{
            \draw(0,-.3)--(0,0)to[out=90,in=90](1,0);
            \draw(0,1)to[out=-90,in=-90](1,1)node[small matrix]{$A$}to[out=90,in=90](2,1)--(2,0)to[out=-90,in=-90](1,0);
            \draw(0,1)node[small matrix]{$A$}--(0,1.4);\node[ciliatednode=170] at(.5,.3){};\node[ciliatednode=190] at(.5,.7){};}
        =0
    \end{equation}
is the equation $A^2-\tr(A)A+\det(A)I=0$. Taking the trace (or connecting both strands instead of just one), one obtains the relationship
    \begin{equation}\label{eq:ch2tr}
        \tr(A^2)-\tr(A)^2+2\det(A)=0.
    \end{equation}
\begin{exercise}
    Use \eqref{eq:ch2tr} to find a formula for $\det(A+B)$.
\end{exercise}

Other versions may be obtained by changing the labels:
    \begin{equation}\label{eq:ch2b-diagram}
        \tikz[xscale=.6]{
            \draw(0,-.3)--(0,0)to[wavyup](1,1)node[small matrix]{$A$}to[out=90,in=90](2,1)--(2,0)to[out=-90,in=-90](1,0);
            \draw(1,0)to[wavyup](0,1)node[small matrix]{$B$}--(0,1.4);}
        -\tikz[xscale=.6]{
            \draw(0,-.3)--(0,0)to[wavyup](.1,.5)to[wavyup](0,1)node[small matrix]{$B$}--(0,1.4);
            \draw(1,0)to[wavyup](.9,.5)to[wavyup](1,1)node[small matrix]{$A$}to[out=90,in=90](2,1)--(2,0)to[out=-90,in=-90](1,0);
            }
        -\tikz[xscale=.6]{
            \draw(0,-.3)--(0,0)to[out=90,in=90](1,0);
            \draw(0,1)to[out=-90,in=-90](1,1)node[small matrix]{$A$}to[out=90,in=90](2,1)--(2,0)to[out=-90,in=-90](1,0);
            \draw(0,1)node[small matrix]{$B$}--(0,1.4);\node[ciliatednode=170] at(.5,.3){};\node[ciliatednode=190] at(.5,.7){};}
        =0
    \end{equation}
is the equation $BA-\tr(A)B+\det(A)BA^{-1}=0$.

We can use \eqref{eq:ch2a-diagram} to express \emph{any} diagram with crossings in terms of diagrams without crossings, which are simply products of traces and determinants. Consider the following family of diagrams with $k$ crossings:
\begin{equation}
    \tikz[heighttwo,xscale=.5]{
        \coordinate(top)at(0,2.5){};\coordinate(topr)at(1,2.5){};\coordinate(bot)at(0,-1.5){};\coordinate(botr)at(1,-1.5){};
        \coordinate(aa)at(0,2){};\coordinate(aar)at(1,2){};\coordinate(bb)at(0,1){};\coordinate(bbr)at(1,1){};
        \coordinate(cc)at(0,-.5){};\coordinate(ccr)at(1,-.5){};
        \draw(bot)to[wavyup](ccr)(botr)to[wavyup](cc)node[small matrix]{$A_k$};
        \draw(bb)node[small matrix]{$A_2$}to[wavyup](aar)(bbr)to[wavyup](aa);
        \node[basiclabel]at(.5,.5){$\vdots$};
        \draw(aa)node[small matrix]{$A_1$}to[wavyup](top)(aar)to[wavyup](topr);
    }
\end{equation}
Applying the binor identity, one obtains a sum of $2^k$ diagrams, from which an identity may be read off.

With two matrices, we have
\begin{equation}\label{eq:2mxid}
    \tikz[heighttwo,xscale=.5]{
        \coordinate(top)at(0,2){};\coordinate(topr)at(1,2){};\coordinate(aa)at(0,1.5){};\coordinate(aar)at(1,1.5){};
        \coordinate(bb)at(0,.5){};\coordinate(bbr)at(1,.5){};\coordinate(cc)at(0,-.5){};\coordinate(ccr)at(1,-.5){};
        \draw(cc)to[wavyup](bbr)(ccr)to[wavyup](bb);
        \draw(bb)node[small matrix]{$B$}to[wavyup](aar)(bbr)to[wavyup](aa);
        \draw(aa)node[small matrix]{$A$}to[wavyup](top)(aar)to[wavyup](topr);
    }
    =\tikz[heighttwo,xscale=.5]{
        \draw(cc)to[wavyup](bb)(ccr)to[wavyup](bbr);
        \draw(bb)node[small matrix]{$B$}to[wavyup](aa)(bbr)to[wavyup](aar);
        \draw(aa)node[small matrix]{$A$}to[wavyup](top)(aar)to[wavyup](topr);
    }
    +\tikz[heighttwo,xscale=.5]{
        \draw(cc)to[wavyup](bb)(ccr)to[wavyup](bbr);
        \draw(bb)node[small matrix]{$B$}to[bend left=90](bbr)(aa)to[bend right=90](aar);
        \draw(aa)node[small matrix]{$A$}to[wavyup](top)(aar)to[wavyup](topr);
        \node[ciliatednode=170] at(.5,.8){};\node[ciliatednode=190] at(.5,1.2){};
    }
    +\tikz[heighttwo,xscale=.5]{
        \draw(cc)to[bend left=90](ccr)(bb)to[bend right=90](bbr);
        \draw(bb)node[small matrix]{$B$}to[wavyup](aa)(bbr)to[wavyup](aar);
        \draw(aa)node[small matrix]{$A$}to[wavyup](top)(aar)to[wavyup](topr);
        \node[ciliatednode=170] at(.5,-.2){};\node[ciliatednode=190] at(.5,.2){};
    }
    +\tikz[heighttwo,xscale=.5]{
        \draw(cc)to[bend left=90](ccr)(bb)to[bend right=90](bbr);
        \draw(bb)node[small matrix]{$B$}to[bend left=90](bbr)(aa)to[bend right=90](aar);
        \draw(aa)node[small matrix]{$A$}to[wavyup](top)(aar)to[wavyup](topr);
        \node[ciliatednode=170] at(.5,.8){};\node[ciliatednode=190] at(.5,1.2){};
        \node[ciliatednode=170] at(.5,-.2){};\node[ciliatednode=190] at(.5,.2){};
    }.
\end{equation}
\begin{exercise}
    Show that the closure (trace) of the diagram \eqref{eq:2mxid} is the tautology
    $$\tr(A)\tr(B)=2\tr(AB)-\tr(BA)-\tr(AB)+\tr(A)\tr(B).$$
\end{exercise}
Higher rank cases are more interesting. With three crossings, the result is the identity
    \begin{multline}\label{eq:tr3}
        \tr(ABC)+\tr(ACB)\\
            =\tr(AB)\tr(C)+\tr(A)\tr(BC)+\tr(B)\tr(CA)-\tr(A)\tr(B)\tr(C).
    \end{multline}
With four crossings, the resulting identity is
    \begin{multline}\label{eq:tr4}
        2\tr(ABCD)\\
            =\tr(A)\tr(BCD)+\tr(B)\tr(ACD)+\tr(C)\tr(ABD)+\tr(D)\tr(ABC) \\
            +\tr(AB)\tr(CD)+\tr(AD)\tr(BC)+\tr(A)\tr(B)\tr(C)\tr(D) \\
            -\tr(A)\tr(B)\tr(CD)-\tr(B)\tr(C)\tr(AD)\\
            -\tr(C)\tr(D)\tr(AB)-\tr(A)\tr(D)\tr(BC).
    \end{multline}
In contrast with the standard proofs, which involve finding and summing several variations of the characteristic equation, these identities arise in a completely straightforward fashion. These results are central to parts of $SL(2,\mathbb{C})$ invariant theory \cite{Gol09}.

And the diagrams suggest further patterns that might be worthy of consideration. One can place matrices at different points within the family, or come up with a completely separate family of diagrams. Another inviting question is looking into how these ideas carry over to higher-dimensional cases.
\fi

\section{Concluding Remarks}

Trace diagrams are particularly suited to recognizing and generalizing patterns in vector and matrix algebra. In many cases, they lead to simple or even ``trivial'' proofs. Moreover, linear algebra is decidedly \emph{non}linear, in the sense that matrix algebra does not naturally play out along a 1-dimensional line of text. One can reasonably argue that the diagrammatic notation for the cross product, the trace, and the determinant is a better approximation of the underlying ``truth'' of these concepts than the standard $\bfu\times\bfv$, $\mathrm{tr}(A)$, and $\det(A)$.

Trace diagrams also have a strong potential to help advance certain areas of mathematics. In invariant theory, they simplify the process of generating and understanding trace relations, which is fundamental for recognizing structures of interest in the theory. In another interesting connection, the four-color theorem may be expressed entirely in terms of these diagrams \cite{KT03}. Independent of these applications, there are many unanswered questions about the diagrams themselves, particularly when matrices are restricted to certain groups.

Trace diagrams are not likely to revolutionize linear algebra, but they provide a refreshing perspective and a new way to understand the subject. And one cannot help but wonder what linear algebra would now look like if trace diagrams had been as easy to typeset as rows and columns of numbers.

\section*{Acknowledgments}
I am grateful to Bill Goldman for introducing me to diagrams, to the University of Maryland for providing plenty of chalkboards to doodle on as a graduate student, and to my research students, Steve Morse, James Lee, and Paul Falcone, for helping the point-of-view presented here to coalesce.

\bibliographystyle{plain}
\bibliography{tracediagrams}

\ifnum0=1
\newpage
\appendix
\section*{Additional Exercises}

\subsection*{3-Diagrams}

\begin{exercise}
What vector identity does this diagram represent?
    \begin{equation}
    \tikz[scale=.6]{
        \node[vertex](node)at(.5,2){};
        \draw[](0,0)node[vector]{$\bf u$}to[out=90,in=-135](node);
        \draw[](1,0)node[vector]{$\bf v$}to[out=90,in=-45](node);
        \draw[](2,0)node[vector]{$\bf w$}to(2,2);\draw(2,2)to[out=90,in=90,looseness=1.5](node);
    }
    =\tikz[scale=.6]{
        \node[vertex](node)at(1.5,2){};
        \draw[](0,0)node[vector]{$\bf u$}to(0,2);\draw(0,2)to[out=90,in=90,looseness=1.5](node);
        \draw[](1,0)node[vector]{$\bf v$}to[out=90,in=-135](node);
        \draw[](2,0)node[vector]{$\bf w$}to[out=90,in=-45](node);
    }
    =\tikz[scale=.6]{
        \node[vertex](node)at(1,2.5){};
        \draw[](0,0)node[vector]{$\bf u$}to[out=90,in=-90](1.5,1.75);\draw(1.5,1.75)to[out=90,in=-45](node);
        \draw[](1,0)node[vector]{$\bf v$}to[out=90,in=-90](2.5,2)to(2.5,2)to[out=90,in=90,looseness=2](node);
        \draw(2,0)node[vector]{$\bf w$}to[out=90,in=-90](.5,1.75);\draw[](.5,1.75)to[out=90,in=-135](node);
    }
    =\tikz[scale=.7]{
        \node[vertex](node)at(1,2){};
        \draw[](0,0)node[vector]{$\bf u$}to[out=90,in=-135](node);
        \draw[](1,0)node[vector]{$\bf v$}to[out=90,in=-90](node);
        \draw[](2,0)node[vector]{$\bf w$}to[out=90,in=-45](node);
    }
    \end{equation}
(The reader is encouraged also to guess the meaning of the fourth term.)
\end{exercise}

\begin{exercise}
    Write down the diagrammatic relation encoding the fact that the cross product $\mathbf{a}\times\mathbf{b}$ is orthogonal to both $\mathbf{a}$ and $\mathbf{b}$. What is the corresponding rule for trace diagrams?
\end{exercise}

\begin{exercise}\label{ex:proj}
    Write down the diagrammatic expression of the \emph{vector projection formula} for projecting a vector $\mathbf{a}$ onto the vector $\mathbf{b}$, $$\mathrm{proj}_{\mathbf{b}} \mathbf{a}=\frac{\mathbf{a}\cdot\mathbf{b}}{|\mathbf{b}|^2}\mathbf{b}.$$
\end{exercise}

\begin{exercise}
    Use diagrams to prove that $|\mathbf{a}\times\mathbf{b}|^2 = |\mathbf{a}|^2|\mathbf{b}|^2 - (\mathbf{a}\cdot\mathbf{b})^2.$
\end{exercise}

\begin{exercise}\label{ex:sincos}
    Suppose the magnitude of a vector is represented with a bar across the end of an open strand, so that
        $$|\mathbf{a}| \equiv \tikz{\draw[](0,0)node[vector]{$\bf a$}to(0,1);\draw[very thick](-.25,1)to(.25,1);}$$
    (a) Express the relationship $|\mathbf{a}|^2=\mathbf{a}\cdot\mathbf{a}$ diagrammatically. (b) Express the cosine angle formula $\cos(\theta)=\frac{\mathbf{a}\cdot\mathbf{b}}{|\mathbf{a}| |\mathbf{b}|}$ diagrammatically (leave the cosine part as is). Do the same for the sine formula $\sin(\theta) = \frac{|\mathbf{a}\times\mathbf{b}|}{|\mathbf{a}| |\mathbf{b}|}$.
\end{exercise}

\begin{exercise}
    Suppose that vector normalization $\mathbf{a} \to \frac{\mathbf{a}}{|\mathbf{a}|}$ is represented diagrammatically by a double hat along the strand, so that
        $$\frac{\mathbf{a}}{|\mathbf{a}|} \equiv \tikz{\draw[](0,0)node[vector]{$\bf a$}to(0,1);\draw[thin](-.15,.6)to(0,.85)to(.15,.6)(-.15,.45)to(0,.7)to(.15,.45);}$$
    Using this notation, find simplified expressions for your answers to exercises \ref{ex:proj} and \ref{ex:sincos}(b).
\end{exercise}

\subsection*{Trace Diagram Computations (Vectors)}

\begin{exercise}
    Use inferred summation to show that $\tikz{\draw(0,.5)circle(.5);} = n$.
\end{exercise}

\begin{exercise}\label{ex:dotprod}
    Prove that $\tikz{\draw(0,.2)node[vector]{$\bfu$}to[out=90,in=90,looseness=2](1,.2)node[vector]{$\bfv$};} = \bfu\cdot\bfv$
    two different ways: first, break up the vectors; and second, use inferred summation on the inside edge.
\end{exercise}

\begin{exercise}\label{ex:det3-b}
    Evaluate the following diagram directly to verify the identity:
    $$
    \tikz[scale=.7]{
        \node[ciliatednode=170](node)at(1,2){};
        \draw[](0,0)node[vector]{$\bf u$}to[out=90,in=-135](node);
        \draw[](1,0)node[vector]{$\bf v$}to[out=90,in=-90](node);
        \draw[](2,0)node[vector]{$\bf w$}to[out=90,in=-45](node);
    }
    =
    \begin{vmatrix}
        u_1 & v_1 & w_1 \\
        u_2 & v_2 & w_2 \\
        u_3 & v_3 & w_3
    \end{vmatrix}
    =(\bfu\times\bfv)\cdot\bfw.
    $$
    Note that this answers a question posed in Exercise \ref{ex:diag-3vec}.
\end{exercise}

\begin{exercise}
    Verify the equality
    $$\tikz{
        \node[vertex]at(0,.5){}edge node[vector,pos=1]{$1$}(0,1)
            edge[bend right]node[vector,pos=1]{$\bfu$}(-.5,0)
            edge[bend left]node[vector,pos=1]{$\bfv$}(.5,0);}
    = u_2 v_3 - u_3 v_2.$$
    Go on to show, using the techniques for evaluating open-strand diagrams, that $\tikz{\node[vertex]at(0,.5){}edge(0,1)edge[bend right]node[vector,pos=1]{$\bfu$}(-.5,0)edge[bend left]node[vector,pos=1]{$\bfv$}(.5,0);} = \bfu \times \bfv$.
\end{exercise}

\begin{exercise}
    Show that
        $\tikz{\draw(0,0)to(0,.3)node[ciliatednode=180]{}to(0,.7)node[ciliatednode=180]{}to(0,1);}
        = -\tikz{\draw(0,0)to(0,1);}$
    and $\tikz{\draw(0,0)to(0,.3)node[ciliatednode=180]{}to(0,.7)node[ciliatednode=0]{}to(0,1);}
        = \tikz{\draw(0,0)to(0,1);}$.
\end{exercise}

\begin{exercise}\label{ex:tensor}
    When a diagram has more than one vector output, the output may be expressed in terms of coefficients of
    each possible labeling of the outputs by basis elements. We can represent these using a \emph{tensor summation}.
    (If you have not seen tensors before, think of them for now as simply placeholders for the output labels.)
    For example, for a diagram with two outputs, the summation
        $a \bse_1\otimes\bse_1 + b\bse_1\otimes\bse_2 + c\bse_2\otimes\bse_1 + d\bse_2\otimes\bse_2$
    indicates that $a$ is the coefficient of the diagram with both output strands labeled 1; $b$ is the coefficient of the diagram
    with first output strand labeled 1 and the second labeled 2; $c$ and $d$ are defined similarly.

    Show that
        $\tikz[xscale=.6]{\draw(0,.5)to[out=-90,in=-90](1,.5);}=\bse_1\otimes\bse_1 + \bse_2\otimes\bse_2$
    and that
        $\tikz[xscale=.6]{\draw(0,.5)to[out=-90,in=-90](1,.5);\node[ciliatednode=190] at(.5,.2){};} = \bse_2\otimes\bse_1 - \bse_1\otimes\bse_2$
    by enumerating all possible labelings of the outputs and computing the coefficient in each case.
\end{exercise}

\begin{exercise}
    (a) Using the previous problem, show that
        $$\tikz[xscale=.6]{\draw(0,0)node[vector,scale=.8]{$\mathbf{u}$}to[out=90,in=90](1,0)node[vector,scale=.8]{$\mathbf{v}$};\draw(0,1)to[out=-90,in=-90](1,1);
            \node[ciliatednode=170] at(.5,.3){};\node[ciliatednode=190] at(.5,.7){};}
        = \det[\mathbf{u}\:\:\mathbf{v}] (\bse_2\otimes\bse_1 - \bse_1\otimes\bse_2).$$
        (Remember that disjoint diagrams are multiplied... see the text for the computation of the bottom diagram.)

    (b) Use your answer to the previous problem to prove the binor identity
    \begin{equation}
    \tikz[xscale=.6]{\draw(0,0)to[wavyup](1,1);\draw(1,0)to[wavyup](0,1);}
    -\tikz[xscale=.6]{\draw(0,0)to[wavyup](.1,.5)to[wavyup](0,1);\draw(1,0)to[wavyup](.9,.5)to[wavyup](1,1);}
    =\tikz[xscale=.6]{\draw(0,0)to[out=90,in=90](1,0);\draw(0,1)to[out=-90,in=-90](1,1);
        \node[ciliatednode=170] at(.5,.3){};\node[ciliatednode=190] at(.5,.7){};}.
    \end{equation}
    Begin by picking specific inputs $\bfu$ and $\bfv$. Then evaluate the left side of the equation as a tensor summation (see Exercise \ref{ex:tensor}).
\end{exercise}

\subsection*{Trace Diagram Computations (Matrices)}

\begin{exercise}
    Let $A=\begin{pmatrix}a&b\\c&d\end{pmatrix}$. (a) Evaluate
        \tikz{\draw(0,.5)circle(.5);\node[small matrix]at(-.5,.5){$A$};\node[ciliatednode=200]at(0,0){};}
    and
        \tikz{\draw(0,.5)circle(.5);\node[small matrix]at(-.5,.5){$A$};\node[small matrix]at(.5,.5){$A$};
            \node[ciliatednode=200]at(0,0){};\node[ciliatednode=160]at(0,1){};}.
    (b) Show that
        $\tikz[xscale=.6]{\draw(-.5,.1)to(-.5,.5)to[out=90,in=90](.5,.5)to(.5,.1);\node[small matrix]at(-.5,.4){$A$};
            \node[small matrix]at(.5,.4){$A$};\node[ciliatednode=160]at(0,.78){};}
        = \det(A) \tikz[xscale=.6]{\draw(-.5,.1)to(-.5,.5)to[out=90,in=90](.5,.5)to(.5,.1);\node[ciliatednode=160]at(0,.78){};}$.
\end{exercise}

\begin{exercise}
    Use dot products (see Exercise \ref{ex:dotprod}) to evaluate
        $\tikz[xscale=.6]{\draw(-.5,.1)to(-.5,.5)to[out=90,in=90](.5,.5)to(.5,.1);\node[small matrix]at(-.5,.4){$A$};}$
    and
        $\tikz[xscale=.6]{\draw(-.5,.1)to(-.5,.5)to[out=90,in=90](.5,.5)to(.5,.1);\node[small matrix]at(.5,.4){$A$};}$.
    Explain why, in general, they are not equal.
\end{exercise}

\begin{exercise}
    (a) Use inferred summation (between the matrix and the input vector) in the diagram \tikz{\draw(0,0)node[vector]{$\bfv$}to[with small matrix={$A$}](0,1);} to verify the matrix-vector multiplication rule in basic linear algebra, i.e. the $i$th coefficient of the
    vector $A\mathbf{v}$ is $\sum_{j=1}^n a_{ij} v_j$. (b) Repeat the process with \tikz{\draw(0,-.1)to[with small matrix={$B$}](0,.5);\draw(0,.5)to[with small matrix={$A$}](0,1.1);} to verify the matrix multiplication rule,
    $(AB)_{ij}=\sum_{k=1}^n A_{ik} B_{kj}$.
\end{exercise}

\begin{exercise}
    (a) Show that
    $$
    \tikz[scale=.7,xscale=.7]{
        \node[ciliatednode=170](node)at(0,2){};
        \draw(-2,0)node[vector,scale=.7]{$\bfu_1$}to[out=90,in=-160](node);
        \draw(-1,0)node[vector,scale=.7]{$\bfu_2$}to[out=90,in=-135](node);
        \draw(-1,.5)to[dotdotdot](2,.5);
        \draw(2,0)node[vector,scale=.7]{$\bfu_n$}to[out=90,in=-20](node);
    }
    = \det [\bfu_1 \: \bfu_2 \: \cdots \: \bfu_n].
    $$

    (b) Using the fact that $\det(AB)=\det(A)\det(B)$, conclude that for an $n\times n$ matrix $A$,
    $$
    \tikz[scale=.7,xscale=.7]{
        \node[ciliatednode=170](node)at(0,2){};
        \draw[with small matrix=$A$](-2,0)node[vector,scale=.7]{$\bfu_1$}to[out=90,in=-160](node);
        \draw[with small matrix=$A$](-1,0)node[vector,scale=.7]{$\bfu_2$}to[out=90,in=-135](node);
        \draw(-1,.5)to[dotdotdot](2,.5);
        \draw[with small matrix=$A$](2,0)node[vector,scale=.7]{$\bfu_n$}to[out=90,in=-20](node);
    }
    = \det(A)
    \tikz[scale=.7,xscale=.7]{
        \node[ciliatednode=170](node)at(0,2){};
        \draw(-2,0)node[vector,scale=.7]{$\bfu_1$}to[out=90,in=-160](node);
        \draw(-1,0)node[vector,scale=.7]{$\bfu_2$}to[out=90,in=-135](node);
        \draw(-1,.5)to[dotdotdot](2,.5);
        \draw(2,0)node[vector,scale=.7]{$\bfu_n$}to[out=90,in=-20](node);
    },
    $$
    and consequently
    $$
    \tikz[scale=.7,xscale=.7]{
        \node[ciliatednode=170](node)at(0,2){};
        \draw[with small matrix=$A$](-2,0)to[out=90,in=-160](node);
        \draw[with small matrix=$A$](-1,0)to[out=90,in=-135](node);
        \draw(-1,.5)to[dotdotdot](2,.5);
        \draw[with small matrix=$A$](2,0)to[out=90,in=-20](node);
    }
    = \det(A)
    \tikz[scale=.7,xscale=.7]{
        \node[ciliatednode=170](node)at(0,2){};
        \draw(-2,0)to[out=90,in=-160](node);\draw(-1,0)to[out=90,in=-135](node);\draw(-1,.5)to[dotdotdot](2,.5);\draw(2,0)to[out=90,in=-20](node);
    }.
    $$

    (c) Finally, show that
    $$
    \tikz[heightoneonehalf]{
            \node[ciliatednode=140](topnode)at(0,1.5){};
            \node[ciliatednode=220](bottomnode)at(0,0){};
            \draw[with small matrix={$A$}](bottomnode)arc(270:90:.75);
            \draw[with small matrix={$A$}](bottomnode)arc(-90:90:.75);
            \draw[with small matrix={$A$}](bottomnode)to[out=135,in=-135](topnode);
            \draw[dotdotdot](-.25,.75)--(.65,.75);
        }
    =\pm n! \det(A).
    $$
    Here, you might proceed by using inferred summation between the bottom node and the matrix labels.
\end{exercise}

\subsection*{Trace Relations}

\begin{exercise}\label{ex:2diagram}
    Redefine local maxima/minima in diagrams by setting
    $\tikz[xscale=.6]{\draw(-.5,.1)to(-.5,.5)to[out=90,in=90](.5,.5)to(.5,.1);}
        \equiv \tikz[xscale=.6]{\draw(-.5,.1)to(-.5,.5)to[out=90,in=90](.5,.5)to(.5,.1);\node[ciliatednode=160]at(0,.78){};}$
    and
    $\tikz[xscale=.6]{\draw(-.5,.9)to(-.5,.5)to[out=-90,in=-90](.5,.5)to(.5,.9);}
        \equiv \tikz[xscale=.6]{\draw(-.5,.9)to(-.5,.5)to[out=-90,in=-90](.5,.5)to(.5,.9);\node[ciliatednode=100]at(0,.22){};}$.
    Use the results from the exercises above to verify the following identities:
    \begin{align}
        \tikz{\draw(0,.5)circle(.5);\node[small matrix]at(-.5,.5){$A$};} &= \tr(A), \\
        \tikz{\draw(0,.5)circle(.5);\node[small matrix]at(-.5,.5){$A$};\node[small matrix]at(.5,.5){$A$};} &= 2\det(A), \\ \tikz[xscale=.6]{\draw(-.5,.1)to(-.5,.5)to[out=90,in=90](.5,.5)to(.5,.1);\node[small matrix]at(-.5,.4){$A$};
            \node[small matrix]at(.5,.4){$A$};}
            &= \det(A) \tikz[xscale=.6]{\draw(-.5,.1)to(-.5,.5)to[out=90,in=90](.5,.5)to(.5,.1);}, \\
        \tikz[xscale=.6]{\draw(0,0)to[wavyup](1,1);\draw(1,0)to[wavyup](0,1);}
            -\tikz[xscale=.6]{\draw(0,0)to[wavyup](.1,.5)to[wavyup](0,1);\draw(1,0)to[wavyup](.9,.5)to[wavyup](1,1);}
            +\tikz[xscale=.6]{\draw(0,0)to[out=90,in=90](1,0);\draw(0,1)to[out=-90,in=-90](1,1);}
            &=0, \\
        \tikz{\draw(0,0)to(0,.6)to[out=90,in=90](.3,.6)to(.3,.4)to[out=-90,in=-90](.6,.4)to(.6,1);}
            &= - \tikz{\draw(0,0)to[out=90,in=-90](.3,1);}.
    \end{align}
\end{exercise}

\begin{exercise}
    Use the identities in the previous problem to describe the following diagram as the product of a single trace and a single determinant:
    $$
    \tikz[heighttwo,xscale=.5]{
        \draw(0,0)to[with small matrix={$C$}](0,1);
        \draw(0,1)to[with small matrix={$D$}](0,2)to[out=90,in=90](1,2);
        \draw(1,2)to[with small matrix={$A$}](1,1);
        \draw(1,1)to[out=-90,in=90](2,0)to[out=-90,in=-90](3,0);
        \draw(3,0)to[with small matrix={$B$}](3,2);
        \draw(3,2)to[out=90,in=90](2,2)node[small matrix]{$B$}to[out=-90,in=90](1,0)to[out=-90,in=-90](0,0);}
    .$$
\end{exercise}

\begin{exercise}
    Show that the closure (trace) of the diagram
    $$
    \tikz[heighttwo,xscale=.5]{
        \coordinate(top)at(0,2){};\coordinate(topr)at(1,2){};\coordinate(aa)at(0,1.5){};\coordinate(aar)at(1,1.5){};
        \coordinate(bb)at(0,.5){};\coordinate(bbr)at(1,.5){};\coordinate(cc)at(0,-.5){};\coordinate(ccr)at(1,-.5){};
        \draw(cc)to[wavyup](bbr)(ccr)to[wavyup](bb);
        \draw(bb)node[small matrix]{$B$}to[wavyup](aar)(bbr)to[wavyup](aa);
        \draw(aa)node[small matrix]{$A$}to[wavyup](top)(aar)to[wavyup](topr);
    }
    =\tikz[heighttwo,xscale=.5]{
        \draw(cc)to[wavyup](bb)(ccr)to[wavyup](bbr);
        \draw(bb)node[small matrix]{$B$}to[wavyup](aa)(bbr)to[wavyup](aar);
        \draw(aa)node[small matrix]{$A$}to[wavyup](top)(aar)to[wavyup](topr);
    }
    +\tikz[heighttwo,xscale=.5]{
        \draw(cc)to[wavyup](bb)(ccr)to[wavyup](bbr);
        \draw(bb)node[small matrix]{$B$}to[bend left=90](bbr)(aa)to[bend right=90](aar);
        \draw(aa)node[small matrix]{$A$}to[wavyup](top)(aar)to[wavyup](topr);
        \node[ciliatednode=170] at(.5,.8){};\node[ciliatednode=190] at(.5,1.2){};
    }
    +\tikz[heighttwo,xscale=.5]{
        \draw(cc)to[bend left=90](ccr)(bb)to[bend right=90](bbr);
        \draw(bb)node[small matrix]{$B$}to[wavyup](aa)(bbr)to[wavyup](aar);
        \draw(aa)node[small matrix]{$A$}to[wavyup](top)(aar)to[wavyup](topr);
        \node[ciliatednode=170] at(.5,-.2){};\node[ciliatednode=190] at(.5,.2){};
    }
    +\tikz[heighttwo,xscale=.5]{
        \draw(cc)to[bend left=90](ccr)(bb)to[bend right=90](bbr);
        \draw(bb)node[small matrix]{$B$}to[bend left=90](bbr)(aa)to[bend right=90](aar);
        \draw(aa)node[small matrix]{$A$}to[wavyup](top)(aar)to[wavyup](topr);
        \node[ciliatednode=170] at(.5,.8){};\node[ciliatednode=190] at(.5,1.2){};
        \node[ciliatednode=170] at(.5,-.2){};\node[ciliatednode=190] at(.5,.2){};
    }.
    $$
    is the expression
    $$\tr(A)\tr(B)=2\tr(AB)-\tr(BA)-\tr(AB)+\tr(A)\tr(B).$$
    Conclude that $\tr(AB)=\tr(BA)$.
\end{exercise}

\begin{exercise}
    Prove \eqref{eq:tr3} by applying the binor identity in Exercise \ref{ex:2diagram} to the crossings in
    $$\tikz[heighttwo,xscale=.5]{
        \draw(0,0)to[out=90,in=-90](1,.67)to[out=90,in=-90](0,1.33)node[small matrix]{$B$}to[out=90,in=-90](1,2)to(1,2.5);
        \draw(1,0)to[out=90,in=-90](0,.67)node[small matrix]{$C$}to[out=90,in=-90](1,1.33)to[out=90,in=-90](0,2)node[small matrix]{$A$}to(0,2.5);
    }.$$ (Assume that the top strands are connected to the bottom strands, so the diagram is closed.)
\end{exercise}

\begin{exercise}
    Prove \eqref{eq:tr4}.
\end{exercise}
\fi

\end{document}